\pdfoutput=1
\documentclass[12pt,a4paper]{scrartcl} 
\usepackage[utf8]{inputenc}
\usepackage[english]{babel}
\usepackage{indentfirst}
\usepackage{misccorr}
\usepackage{graphicx}
\usepackage{amsmath}
\usepackage[usenames]{color}
\begin{document}
	\begin{center}
		\large
		\textbf{The real parts of the nontrivial Riemann zeta function zeros}
		\\
		Igor Turkanov
		\\
		to my love and wife Mary 
	\end{center}
	\smallskip
	\begin{quote}
		\large
		ABSTRACT
		\\	
		\\
		This theorem is based on holomorphy of studied functions and the fact that near a singularity point the real part of some rational function can take an arbitrary preassigned value.
	\end{quote}
	\begin{quote}
		\normalsize
		The colored markers are as follows:
		\\
		$\textcolor{red}{\bullet}$ -  assumption or a fact, which is not proven at present;\\
		$\textcolor{yellow}{\bullet}$ - the statement, which requires additional attention;\\
		$\textcolor{green}{\bullet}$ - statement, which is proved earlier or clearly undestandable.
		\\
	\end{quote}
	\large
	THEOREM
	\\
	\\
	\marginparsep = 10pt
	\marginparwidth = 10pt 
	\reversemarginpar
	\marginpar{$\textcolor{red}{\bullet}$}
	The real parts of all the nontrivial Riemann zeta function zeros $ \rho $ are equal $ Re\left(\rho\right)=\dfrac{1}{2}. $
	\\
	\\
	\\
	PROOF:
	\\
	\\\marginpar{$\textcolor{green}{\bullet}$}
	In relation to $ \zeta\left(s\right) $ - Zeta function of Riemann is known~\cite[p.~5]{Titchmarsh-1987} two equations each of which can serve as its definition:
	\begin{equation} \label{function_eq_Eiler}
	\centering
	\zeta \left(s\right)=\sum_{n=1}^{\infty} \dfrac{1}{n^s},\qquad\zeta \left(s\right)=\prod_{n=1}^{\infty} \left( 1-\dfrac{1}{p_n^s}\right)^{-1} ,\qquad Re\left(s\right) >1,
	\end{equation}
	where $ p_1, p_2, \ldots, p_n, \ldots $ is a series of primes.
	\\
	\\\marginpar{$\textcolor{green}{\bullet}$}
	According to the functional equality~\cite[p.~22]{Titchmarsh-1987}, ~\cite[p.~8-11]{Karatsuba-1992} by part $\Gamma\left(s\right)$ is the Gamma function:
	\begin{equation} \label{function_eq}
	\centering
	\Gamma \left(\dfrac{s}{2}\right){\pi }^{-\dfrac{s}{2}}\zeta \left(s\right)=\Gamma \left(\dfrac{1-s}{2}\right){\pi }^{-\dfrac{1-s}{2}}\zeta \left(1-s\right),\qquad Re\left(s\right) >0.
	\end{equation}
	\marginpar{$\textcolor{green}{\bullet}$}
	From~\cite[p.~8-11]{Karatsuba-1992}\; $\zeta \left({\bar{s}}\right)=\overline{\zeta\left(s\right)}$, it means that $\forall \rho =\sigma +it$: $\zeta\left(\rho\right)=0$ and $0\leqslant \sigma \leqslant 1$ we have:
	\begin{equation}\label{eq_zero}
	\zeta \left({\bar{\rho}}\right)=\zeta \left({1-\rho}\right)=\zeta \left({1-\bar{\rho}}\right)=0
	\end{equation}
	\marginpar{$\textcolor{green}{\bullet}$}
	From~\cite{Valle-Poussin-1897},~\cite[p.~128]{Smith-1994},~\cite[p.~45]{Titchmarsh-1987} \; we know that  $\zeta \left({s}\right)$ has no nontrivial zeros on the line $\sigma = 1$ and consequently on the line $\sigma = 0$ also, in accordance with~\eqref{eq_zero} they don't exist.
	\\
	\\
	\marginpar{$\textcolor{yellow}{\bullet}$}
	Let's denote the set of nontrivial zeros  $\zeta \left({s}\right)$ through $ \mathcal{P} $ (multiset with consideration of multiplicitiy):
	\[ \mathcal{P} \stackrel{\scriptscriptstyle\mathrm{def}}{=} \left\lbrace \rho:\; \zeta\left(\rho\right)=0,\; \rho =\sigma +it,\; 0<\sigma <1  \right\rbrace. \]
	\begin{eqnarray}
	\mbox{And:}\;\mathcal{P}_1 \stackrel{\scriptscriptstyle\mathrm{def}}{=} \left\lbrace \rho:\; \zeta\left(\rho\right)=0,\; \rho =\sigma +it,\; 0<\sigma <\dfrac{1}{2}  \right\rbrace,\nonumber\\
	\mathcal{P}_2 \stackrel{\scriptscriptstyle\mathrm{def}}{=} \left\lbrace \rho:\; \zeta\left(\rho\right)=0,\; \rho =\dfrac{1}{2} +it \right\rbrace\,,\quad\qquad\qquad\nonumber\\
	\mathcal{P}_3 \stackrel{\scriptscriptstyle\mathrm{def}}{=} \left\lbrace \rho:\; \zeta\left(\rho\right)=0,\; \rho =\sigma +it,\; \dfrac{1}{2}<\sigma < 1  \right\rbrace\nonumber.
	\end{eqnarray}	
	Then: \[ \mathcal{P}= \mathcal{P}_1 \cup  \mathcal{P}_2 \cup  \mathcal{P}_3\;\;\; \mbox{and}\;\; \mathcal{P}_1 \cap  \mathcal{P}_2 =  \mathcal{P}_2 \cap  \mathcal{P}_3 =  \mathcal{P}_1 \cap  \mathcal{P}_3 = \varnothing , \]
	\[  \mathcal{P}_1 =\varnothing \Leftrightarrow  \mathcal{P}_3=\varnothing.\]
	\marginpar{$\textcolor{green}{\bullet}$}
	Hadamard's theorem (Weierstrass preparation theorem) about the\\ decomposition of function through the roots gives us the following result\\~\cite[p.~30]{Titchmarsh-1987}, ~\cite[p.~31]{Karatsuba-1992},~\cite{Voros-1987}:\\
	\begin{equation}\label{eq_Adamar}
	\zeta \left({s}\right)= \dfrac{{\pi }^{\frac{s}{2}}{e}^{as}}{s\left(s-1\right)\Gamma \left(\frac{s}{2}\right)}\prod _{\rho \in \mathcal{P}}\left(1-\frac{s}{\rho }\right){e}^{\frac{s}{\rho }},\qquad Re\left(s\right) >0 
	\end{equation}
	\[ a=\mathrm{ln}2\sqrt{\pi }-\dfrac{\gamma }{2}-1,\; \gamma - \mbox{Euler's constant and} \]	
	\begin{equation}\label{eq_PrimeLn}
	\dfrac{{\zeta }^{\prime }\left(s\right)}{\zeta \left(s\right)}=\dfrac{1}{2}\mathrm{ln}\pi + a - \dfrac{1}{s} + \dfrac{1}{1-s} -\dfrac{1}{2}\dfrac{{\Gamma }^{\prime }\left(\dfrac{s}{2}\right)}{\Gamma \left(\dfrac{s}{2}\right)}+\sum _{\rho \in \mathcal{P}}\left(\dfrac{1}{s-\rho }+\dfrac{1}{\rho }\right)
	\end{equation}
	\\
	\marginpar{$\textcolor{green}{\bullet}$}
	According to the fact that  $\;\;\dfrac{{\Gamma }^{\prime }\left(\dfrac{s}{2}\right)}{\Gamma \left(\dfrac{s}{2}\right)}\;\; $   - Digamma function of~\cite[p.~31]{Titchmarsh-1987},\\ ~\cite[p.~23]{Karatsuba-1992} we have:
	\begin{equation}\label{eq_PrimeLnDigamma}
	\dfrac{{\zeta }^{\prime }\left(s\right)}{\zeta \left(s\right)}= \dfrac{1}{1-s} +\sum _{\rho \in \mathcal{P}}\left(\dfrac{1}{s-\rho }+\dfrac{1}{\rho }\right)+\sum _{n=1}^{\infty }\left(\dfrac{1}{s+2n}-\dfrac{1}{2n}\right)+C,
	\end{equation}
	\[ C=const. \]
	\\
	\marginpar{$\textcolor{green}{\bullet}$}
	From~\cite[p.~160]{Edwards-1974}, ~\cite[p.~272]{Lehmer-1988}, ~\cite[p.~81]{Davenport-1980}:
	\begin{equation}\label{eq_RHO}
	\qquad \sum _{\rho \in \mathcal{P}}\dfrac{1}{\rho }=1+\dfrac{\gamma }{2}-\mathrm{ln}2\sqrt{\pi }=0,0230957\dotsc 
	\end{equation}
	\\
	\marginpar{$\textcolor{yellow}{\bullet}$}
	Indeed, from~\eqref{eq_zero}:
	\[  \sum _{\rho \in \mathcal{P}}\dfrac{1}{\rho }=\dfrac{1}{2}\sum _{\rho \in \mathcal{P}}\left(\dfrac{1}{1-\rho }+\dfrac{1}{\rho }\right). \]
	\\
	\marginpar{$\textcolor{green}{\bullet}$}
	From~\eqref{eq_PrimeLn}:
	\[2 \sum _{\rho \in \mathcal{P}}\dfrac{1}{\rho }=\underset{s\to 1}{lim}\left(\dfrac{{\zeta }^{\prime }\left(s\right)}{\zeta \left(s\right)}-\dfrac{1}{1-s}+\dfrac{1}{s}-a-\dfrac{1}{2}\mathrm{ln}\pi+\dfrac{1}{2}\dfrac{{\Gamma }^{\prime }\left(\dfrac{s}{2}\right)}{\Gamma \left(\dfrac{s}{2}\right)}\right). \]
	\\
	\marginpar{$\textcolor{green}{\bullet}$}
	Also it's known, for example, from ~\cite[p.~49]{Titchmarsh-1987}, ~\cite[p.~98]{Davenport-1980} that the number of nontrivial zeros of $\rho =\sigma +it$ in strip $ 0<\sigma <1 $, the imaginary parts of which $t$ are less than some number $T>0$ is limited, i.e.,
	\begin{equation*}\label{eq_SetZeroes}
	\| \left\lbrace \rho:\; \rho \in \mathcal{P},\; \rho =\sigma +it,\; \left| t\right|  <T  \right\rbrace\|<\infty.
	\end{equation*}
	\\
	\marginpar{$\textcolor{yellow}{\bullet}$}
	Indeed, it can be presented that on the contrary the sum of  $  \sum _{\rho \in \mathcal{P}}\dfrac{1}{\rho } $ would have been unlimited.
	\\
	\pagebreak
	\\
	\marginpar{$\textcolor{yellow}{\bullet}$}
	Thus  $\forall\; T>0 \; \exists\; \delta_x>0, \; \delta_y>0$ such that
	in area $ 0<t\leqslant \delta_y,\; 0<\sigma\leqslant \delta_x $ there are no zeros $ \rho =\sigma +it \in \mathcal{P} $.
	\\
	\\
	Let's consider random root $ q \in \mathcal{P}$.
	\\
	\\
	Let's denote $k(q)$ the multiplicity of the root $q$.
	\\
	\\
	Let's examine the area $Q\left(R\right)\stackrel{\scriptscriptstyle\mathrm{def}}{=} \left\{s:\left\| s-q\right\|  \leqslant R, R>0\right\}$.
	\\
	\\
	\marginpar{$\textcolor{yellow}{\bullet}$}
	From the fact of finiteness of set of nontrivial zeros $ \zeta(s) $ in the limited area follows $\exists\; R>0$, such that $Q(R)$ does not contain any root from $ \mathcal{P}$ except $q$ and also does not intersect with the axes of coordinates.
	\\
	\renewcommand{\figurename}{Fig.}
	\begin{figure}[h]	
		\begin{picture}(330,100)(-100,0)
		\put(2,2){\begin{picture}(270,240)%
			\put(15,0){\vector(0,1){110}}
			\put(0,15){\vector(1,0){260}}
			\put(230,2){$Re(s)$} \put(-24,100){$Im(s)$}
			\put(-4,-4){$0$}
			\multiput(30,0)(0,3){37}%
			{\circle*{1}}
			\multiput(220,0)(0,3){37}%
			{\circle*{1}}
			\multiput(120,0)(0,3){37}%
			{\circle*{1}}
			\multiput(0,30)(3,0){86}%
			{\circle*{1}}
			\put(17,-4){$\delta_x$}	
			\put(110,-4){$\frac{1}{2}$}
			\put(210,-4){$1$}
			\put(-12,20){$\delta_y$}
			\put(80,80){\circle{70}}
			\put(80,80){\circle*{2}}
			\put(80,80){\vector(3,1){19}}	\put(81,85){$R$}
			\put(72,81){$q$}
			\put(50,40){\circle*{2}}	
			\put(190,40){\circle*{2}}
			\put(160,80){\circle*{2}}
			\put(120,55){\circle*{2}}	
			\put(120,90){\circle*{2}}
			\end{picture}}
	\end{picture}
	\caption{}
\end{figure}
\\
\marginpar{$\textcolor{green}{\bullet}$}
From~\cite{Abramowitz-1972}, ~\cite[p.~31]{Titchmarsh-1987},~\cite[p.~23]{Karatsuba-1992} we know that the Digamma function $\dfrac{{\Gamma}^{\prime }\left(\dfrac{s}{2}\right)}{\Gamma \left(\dfrac{s}{2}\right)} $ in the area $Q(R)$ has no poles, i.e.,  $\forall s \in Q(R)$
\[\left\| \dfrac{{\Gamma}^{\prime }\left(\dfrac{s}{2}\right)}{\Gamma \left(\dfrac{s}{2}\right)}\right\|  <\infty. \]
Let's denote:
\[	I_{\mathcal{P} }(s)\stackrel{\scriptscriptstyle\mathrm{def}}{=} - \dfrac{1}{s} + \dfrac{1}{1-s}+\sum _{\rho \in \mathcal{P}}\dfrac{1}{s-\rho }  \]
and
\begin{eqnarray} 
I_{\mathcal{P}\setminus \left\lbrace q\right\rbrace }(s)\stackrel{\scriptscriptstyle\mathrm{def}}{=} - \dfrac{1}{s} + \dfrac{1}{1-s}+\sum _{\rho \in \mathcal{P}\setminus \left\lbrace q\right\rbrace}\dfrac{1}{s-\rho }.\nonumber
\end{eqnarray}
\\
\marginpar{$\textcolor{yellow}{\bullet}$}
Hereinafter $\mathcal{P}\setminus \left\lbrace q\right\rbrace \stackrel{\scriptscriptstyle\mathrm{def}}{=} \mathcal{P}\setminus \left\lbrace (q, k(q))\right\rbrace $ (the difference in the multiset).
\\
\\
Also we shall consider the summation - $\sum _{\rho \in \mathcal{P}}\dfrac{1}{s-\rho }  $ and $ \sum _{\rho \in \mathcal{P}\setminus \left\lbrace q\right\rbrace}\dfrac{1}{s-\rho }  $ further as the sum of pairs $\left( \dfrac{1}{s-\rho }+\dfrac{1}{s-(1-\rho) } \right)  $ and $ \sum _{\rho \in \mathcal{P}}\dfrac{1}{\rho } $ as the sum of pairs $\left( \dfrac{1}{\rho }+\dfrac{1}{1-\rho} \right)  $ as a consequence of division of the sum from~\eqref{eq_PrimeLnDigamma} $ \sum _{\rho \in \mathcal{P}}\left(\dfrac{1}{s-\rho }+\dfrac{1}{\rho }\right) $ into $ \sum _{\rho \in \mathcal{P}}\dfrac{1}{s-\rho }+\sum _{\rho \in \mathcal{P}}\dfrac{1}{\rho } $. As specifed in ~\cite{Edwards-1974}, ~\cite{Keiper-1992}, ~\cite{Lehmer-1988}, ~\cite{Titchmarsh-1987}.
\\
\\
\\
\marginpar{$\textcolor{yellow}{\bullet}$}
Let's note that $ I_{\mathcal{P}\setminus \left\lbrace q\right\rbrace }(s) $ is holomorphic function $\forall\; s \in Q(R) $.
\\
\\
Then from~\eqref{eq_PrimeLn} we have:
\begin{eqnarray} 
\dfrac{{\zeta }^{\prime }\left(s\right)}{\zeta \left(s\right)}=\dfrac{1}{2}\mathrm{ln}\pi + a -\dfrac{1}{2}\dfrac{{\Gamma }^{\prime }\left(\dfrac{s}{2}\right)}{\Gamma \left(\dfrac{s}{2}\right)}+\sum _{\rho \in \mathcal{P}}\dfrac{1}{\rho }+I_{\mathcal{P} }(s).\nonumber
\end{eqnarray}
\\
\marginpar{$\textcolor{yellow}{\bullet}$}
And in view of \eqref{eq_Adamar}, \eqref{eq_RHO}:
\begin{equation}\label{eq_Dzeta_Re}
Re\dfrac{{\zeta }^{\prime }\left(s\right)}{\zeta \left(s\right)}=\dfrac{1}{2}\mathrm{ln}\pi +Re\left( -\dfrac{1}{2}\dfrac{{\Gamma }^{\prime }\left(\dfrac{s}{2}\right)}{\Gamma \left(\dfrac{s}{2}\right)}+I_{\mathcal{P} }(s)\right) .
\end{equation}
\\
\marginpar{$\textcolor{yellow}{\bullet}$}
Let's note that from the equality of 
\begin{equation}\label{eq_Minus_rho}
\sum _{\rho \in \mathcal{P}}\dfrac{1}{1-s-\rho }= -\sum _{(1-\rho) \in \mathcal{P}}\dfrac{1}{s-(1-\rho) } = -\sum _{\rho \in \mathcal{P}}\dfrac{1}{s-\rho }
\end{equation}
\\
\marginpar{$\textcolor{yellow}{\bullet}$}
follows that:
\[	I_{\mathcal{P} }(1-s)=-	I_{\mathcal{P} }(s), \; I_{\mathcal{P}\setminus \left\lbrace1- q\right\rbrace }(1-s)= -I_{\mathcal{P}\setminus \left\lbrace q\right\rbrace }(s), \;Re\left(s\right) >0. \]
\\
\\
\marginpar{$\textcolor{green}{\bullet}$}
Besides  \[ I_{\mathcal{P}\setminus \left\lbrace q\right\rbrace }(s)=I_{\mathcal{P} }(s)-\dfrac{k(q)}{s-q} \]
and $ I_{\mathcal{P}\setminus \left\lbrace q\right\rbrace }(s)  $ is limited in the area of $ s \in Q(R) $ as a result of absence of its poles in this area as well as its differentiability in each point of this area.
\\
\\
\marginpar{$\textcolor{green}{\bullet}$}
If in \eqref{eq_PrimeLn} we replace $ s $ with $1-s$ that in view of \eqref{eq_RHO}, in a similar way if we take derivative of the principal logarithm \eqref {function_eq}:
\begin{equation}\label{eq_SumDig}
\dfrac{{\zeta }^{\prime }\left(s\right)}{\zeta \left(s\right)}+\dfrac{{\zeta }^{\prime }\left(1-s\right)}{\zeta \left(1-s\right)}=-\dfrac{1}{2}\dfrac{{\Gamma }^{\prime }\left(\dfrac{s}{2}\right)}{\Gamma \left(\dfrac{s}{2}\right)}-\dfrac{1}{2}\dfrac{{\Gamma }^{\prime }\left(\dfrac{1-s}{2}\right)}{\Gamma \left(\dfrac{1-s}{2}\right)}+\ln\pi,\;Re\left(s\right) >0. 	
\end{equation}
\\
\marginpar{$\textcolor{green}{\bullet}$}
Let's examine a circle with the center in a point  $q$ and radius $r \leqslant R$, laying in the area of $Q(R)$:
\\
\begin{figure}[h]
	\begin{picture}(330,100)(-100,0)
	\put(2,2){\begin{picture}(270,250)%
		\put(15,0){\vector(0,1){110}}
		\put(0,15){\vector(1,0){260}}
		\put(230,2){$Re(s)$} \put(-24,100){$Im(s)$}
		\put(-4,-4){$0$}
		\multiput(30,0)(0,3){37}%
		{\circle*{1}}
		\multiput(220,0)(0,3){37}%
		{\circle*{1}}
		\multiput(120,0)(0,3){37}%
		{\circle*{1}}
		\multiput(0,30)(3,0){86}%
		{\circle*{1}}
		\put(17,-4){$\delta_x$}	
		\put(110,-4){$\frac{1}{2}$}
		\put(210,-4){$1$}
		\put(-9,19){$\delta_y$}
		\put(80,80){\circle{70}}
		\put(80,80){\circle*{2}}
		\put(80,80){\vector(3,1){19}}	
		\put(80,84){$R$}
		\put(74,84){$q$}
		\put(50,40){\circle*{2}}	
		\put(190,40){\circle*{2}}
		\put(160,80){\circle*{2}}
		\put(120,55){\circle*{2}}	
		\put(120,90){\circle*{2}}
		\put(80,80){\circle{33}}
		\put(80,80){\vector(3,-4){10}}
		\put(87,72){$r$}
		\multiput(80,12)(0,4){18}%
		{\circle*{1}}
		\multiput(12,80)(4,0){18}%
		{\circle*{1}}
		\multiput(12,65)(4,0){17}%
		{\circle*{1}}
		\put(75,65){\circle*{2}}  		
		\multiput(75,12)(0,4){14}%
		{\circle*{1}}
		\put(58,55){$m_r$}	
		\put(62,0){${x}_{m_r}$}
		\put(80,0){$\sigma_q$}
		\put(-3,80){$t_q$}	
		\put(-3,58){$y_{m_r}$}
		\end{picture}}
\end{picture}
\caption{}
\end{figure}
\\
\marginpar{$\textcolor{green}{\bullet}$}
For $s=x+iy, \; q=\sigma_q+it_q$
\[ Re\dfrac{k(q)}{s-q}= Re\dfrac{k(q)}{x+iy-\sigma_q-it_q}=\dfrac{k(q)(x-\sigma_q)}{(x-\sigma_q)^2+(y-t_q)^2}=k(q)\dfrac{x-\sigma_q}{r^2}. \]
\\
\\
Let's prove a series of statements:
\\
\pagebreak
\\
\\
\marginpar{$\textcolor{red}{\bullet}$}
STATEMENT A
\\
\\
For an arbitrary nontrivial zero $ q \in \mathcal{P}$ the equality must hold:
\[\forall\; q=\sigma_q+it_q \in \mathcal{P}:  \]
\begin{eqnarray} 
&\dfrac{\partial}{\partial x}\left.Re\left( -\dfrac{1}{2}\dfrac{{\Gamma }^{\prime }\left(\dfrac{x+it_q}{2}\right)}{\Gamma \left(\dfrac{x+it_q}{2} \right)} -\dfrac{1}{2}\dfrac{{\Gamma }^{\prime }\left(\dfrac{1-x-it_q}{2}\right)}{\Gamma \left(\dfrac{1-x-it_q}{2} \right)}\right)\right|_{x=\sigma_q}=0.\nonumber
\end{eqnarray}
\\
PROOF:
\\
\\
\marginpar{$\textcolor{yellow}{\bullet}$}
Let us denote for $ s=x+iy \in Q(R) $:
\begin{align} \label{eq_Def_Tqs}
T_q(s)\stackrel{\scriptscriptstyle\mathrm{def}}{=}\dfrac{{\zeta }^{\prime }\left(s\right)}{\zeta \left(s\right)}-\dfrac{k(q)}{s-q} =\dfrac{1}{2}\mathrm{ln}\pi-\dfrac{1}{2}\dfrac{{\Gamma }^{\prime }\left(\dfrac{s}{2}\right)}{\Gamma \left(\dfrac{s}{2}\right)}+I_{\mathcal{P}\setminus \left\lbrace q\right\rbrace }(s). 
\end{align}
\\
\marginpar{$\textcolor{yellow}{\bullet}$}
The function $ T_q(s) $ is analytic and bounded in the entire neighborhood of $Q(R) $.
\\
\\
\marginpar{$\textcolor{red}{\bullet}$}
Assume that:
\begin{align} \label{eq_NotZero_q}
\dfrac{\partial}{\partial x} Re T_q(x+it_q)\left|_{ x=\sigma_q }\right. \ne 0.
\end{align}
\\
\\
\marginpar{$\textcolor{green}{\bullet}$}
By virtue of the continuity of the function $ \dfrac{\partial}{\partial x} Re T_q(x+iy) $ there exists a neighborhood $ U(q) $ of the point $ q $ such that this function will not turn to $ 0 $, or rather will have the same sign as $ \dfrac{\partial}{\partial x} Re T_q(x+it_q)\left|_{ x=\sigma_q }\right. $ for all points $ s=x+iy \in U(q) $:
\begin{align} \label{eq_NotZero_Uq}
sign \left( \dfrac{\partial}{\partial x} Re T_q(s)\right)=sign\left( \dfrac{\partial}{\partial x} Re T_q(x+it_q)\left|_{ x=\sigma_q }\right.\right)\ne 0.
\end{align}
\\
\\
\marginpar{$\textcolor{green}{\bullet}$}
Hence there must exist a radius $ R_U>0 $ such that the neighborhood $ Q(R_U) $ centered at point $ q $ and radius $ R_U \leqslant R $ is entirely included in the neighborhood $ U(q) $:
\[ \exists\; 0<R_U\leqslant R:\;\; Q(R_U) \subset  U(q). \]
\\
\marginpar{$\textcolor{yellow}{\bullet}$}
That is, for any $ s=x+iy \in Q(R_U) $, the equality \eqref{eq_NotZero_Uq} is satisfied.
\\
\\
\marginpar{$\textcolor{yellow}{\bullet}$}
For arbitrary $ 0 < r \leqslant R_U $ and $ \Delta_r,\;\theta \in \mathbb{R}:\;\; \Delta_r>0,\;\;\theta>0 $ we define the points:
\[a_r \stackrel{\scriptscriptstyle\mathrm{def}}{=}\sigma_q+\Delta_r+i\left( t_q+r\right),  \]
\[v_r \stackrel{\scriptscriptstyle\mathrm{def}}{=}\sigma_q+\Delta_r+\theta\Delta_r +i\left( t_q+r\right).  \]
\\
Here, the points $ a_r $ and $ v_r $ lie on a line parallel to the abscissa axis and passing through the upper point of the circle centered at $ q $ and radius $ r $. The point $ v_r $ is displaced from $ a_r $ along the horizontal axis by the value $ \theta\Delta_r $.
\\
\\
\marginpar{$\textcolor{yellow}{\bullet}$}
Let's find such values of $\Delta_r $ and $\theta $, if they exist, that the equality is satisfied:
\begin{align} \label{eq_Eqviv_ar_vr}
 Re\dfrac{{\zeta }^{\prime }\left(a_r\right)}{\zeta \left(a_r\right)}=Re\dfrac{{\zeta }^{\prime }\left(v_r\right)}{\zeta \left(v_r\right)}.
\end{align}
\\
\marginpar{$\textcolor{green}{\bullet}$}
Which, given the notation \eqref{eq_Def_Tqs}, looks like this:
\[ ReT_q(a_r)+Re\dfrac{k(q)}{a_r-q}= ReT_q(v_r)+Re\dfrac{k(q)}{v_r-q}. \]
\\
\marginpar{$\textcolor{green}{\bullet}$}
And since the values of $ \theta $ are assumed to be positive, based on the mean value theorem (or Lagrange theorem), there must exist a point $ u_r=x_{u_r}+i(t_q+r) $ on the interval $ (a_r,\; v_r) $ such that:
\[ ReT_q(v_r)-ReT_q(a_r)=\theta\Delta_r\dfrac{\partial}{\partial x} Re T_q(x+i(t_q+r))\left|_{ x=x_{u_r} }\right. \]
\pagebreak
\\
\marginpar{$\textcolor{green}{\bullet}$}
Let us show the location of the points under consideration in the following figure:
\\
\begin{figure}[h]
	\begin{picture}(340,350)(-60,0)
	\put(2,2){\begin{picture}(340,350)%
		\put(15,0){\vector(0,1){350}}
		\put(0,15){\vector(1,0){370}}
		\put(358,2){$x$} 
		\put(0,340){$y$}
		\put(180,180){\circle*{3}}
		\put(188,170){$q$}
		\multiput(12,180)(4,0){85}
		{\circle*{1}}
		\multiput(180,8)(0,4){43}
		{\circle*{1}}
		\multiput(180,180)(0,4){38}
		{\circle*{1}}
		\put(177,-5){$\sigma_q$}
		\put(0,178){$t_q$}
		\put(180,180){\vector(-4,1){137}}
		\qbezier(46,227)(23,145)(81,81)
		\qbezier(81,81)(125,39)(180,39)
		\qbezier(46,227)(59,262)(79,281)
		\qbezier(281,81)(235,39)(180,39)
		\qbezier(79,281)(123,321)(180,321)
		\qbezier(314,227)(337,145)(281,81)
		\qbezier(314,227)(301,262)(281,281)
		\qbezier(281,281)(240,321)(180,321)
		\put(180,180){\vector(1,1){55}}
		\multiput(235,235)(3,3){10}
		{\circle*{1}}
		\qbezier(192,192)(199,188)(198,180)
		\put(199,187){$45^{o}$}
		\qbezier(125,125)(153,101)(180,102)
		\qbezier(125,125)(102,151)(102,180)	
		\qbezier(235,125)(207,101)(180,102)
		\qbezier(235,125)(258,151)(258,180)	
		\qbezier(235,235)(207,259)(180,258)
		\qbezier(235,235)(258,209)(258,180)	
		\qbezier(125,235)(153,259)(180,258)
		\qbezier(125,235)(102,209)(102,180)	
		\put(48,200){$R_0$}
		\put(227,218){$r$}
		\multiput(12,258)(5,0){68}
		{\circle*{1}}
		\multiput(220,258)(1,0){65}
		{\circle*{2}}
		\put(-20,255){$t_q+r$}
		\multiput(220,258)(0,-5){51}
		{\circle*{1}}
		\put(202,-5){$\sigma_q+\Delta_r$}
		\multiput(285,258)(0,-5){51}
		{\circle*{1}}
		\put(267,-5){$\sigma_q+\Delta_r+\theta\Delta_r$}	
		\put(285,258){\circle*{4}}
		\put(220,258){\circle*{4}}
		\put(270,258){\circle*{4}}
		\put(252,258){\circle*{4}}
		\put(215,263){$a_r$}
		\put(282,263){$v_r$}
		\put(263,263){$u_{v_r}$}
		\put(246,263){$z_{v_r}$}
		\put(180,180){\vector(-1,2){68}}
		\qbezier(73,73)(120,29)(180,29)
		\qbezier(73,73)(29,120)(29,180)
		\qbezier(287,73)(240,29)(180,29)
		\qbezier(287,73)(331,120)(331,180)
		\qbezier(288,288)(331,245)(331,180)
		\qbezier(288,288)(245,330)(180,330)
		\qbezier(73,290)(115,330)(180,330)
		\qbezier(73,290)(29,247)(29,180)
		\put(105,290){$R_w$}		
		\end{picture}}
\end{picture}
\caption{}
\end{figure}
\\
\\
\marginpar{$\textcolor{green}{\bullet}$}
The inequality \eqref{eq_Eqviv_ar_vr} can be rewritten as follows:
\[   k(q)\dfrac{\Delta_r}{r^2+\Delta_r^2}-k(q)\dfrac{\Delta_r(1+\theta)}{r^2+\left( \Delta_r(1+\theta) \right)^2 }=Re T_q(v_r)-Re T_q(a_r).\]
\\
\marginpar{$\textcolor{yellow}{\bullet}$}
Let's denote:
\[ \alpha_r\stackrel{\scriptscriptstyle\mathrm{def}}{=}\dfrac{1}{k(q)}\dfrac{\partial}{\partial x} Re T_q(x+i(t_q+r))\left|_{ x=x_{u_r} }\right.. \]
\\
\\
\marginpar{$\textcolor{green}{\bullet}$}
Then the equality \eqref{eq_Eqviv_ar_vr} looks like this:
\[ \dfrac{\Delta_r}{r^2+\Delta_r^2}-\dfrac{\Delta_r(1+\theta)}{r^2+\left( \Delta_r(1+\theta) \right)^2 }=\theta\Delta_r\alpha_r. \] 
\\
\marginpar{$\textcolor{green}{\bullet}$}
Since we are looking for $\Delta_r>0 $, the value of $\Delta_r $ we are looking for must satisfy the equation:
\begin{align} 
 \dfrac{1}{r^2+\Delta_r^2}-\dfrac{1+\theta}{r^2+\left( \Delta_r(1+\theta) \right)^2 }=\theta\alpha_r.\nonumber 
\end{align}
\\
\marginpar{$\textcolor{green}{\bullet}$}
Or:
\[ \theta(1+\theta)\Delta_r^2-\theta r^2=\theta\alpha_r\left(r^2+\left( \Delta_r(1+\theta) \right)^2  \right) \left(r^2+\Delta_r^2 \right).  \]
And given the constraint $ \theta>0 $:
\[ (1+\theta)\Delta_r^2-r^2 =\alpha_r\left(r^2+\left( \Delta_r(1+\theta) \right)^2  \right) \left(r^2+\Delta_r^2 \right).  \]
And after opening the brackets:
\begin{align} \label{eq_Ur_t_Dr_ar_fin}
 \alpha_r(1+\theta)^2\Delta_r^4-\left(1+\theta-\alpha_r r^2\left(1+(1+\theta)^2 \right)  \right)\Delta_r^2+r^2+ \alpha_r r^4=0.
\end{align}
\\
\marginpar{$\textcolor{green}{\bullet}$}
Assuming \eqref{eq_NotZero_q} and that for any $ s=x+iy \in Q(R_U)$ the expression \eqref{eq_NotZero_Uq} is satisfied, and given the choice of $ \theta > 0 $ we have a non-zero coefficient at higher degree for any $ 0 < r \leqslant R_U$:
\[ \alpha_r(1+\theta)^2 \ne 0, \]
i.e., the expression \eqref{eq_Ur_t_Dr_ar_fin} is a quadratic equation with respect to the variable $ \Delta_r^2 $, where the discriminant is equal:
\\
\[ D_{\Delta_r^2} \stackrel{\scriptscriptstyle\mathrm{def}}{=}\left( 1+\theta-\alpha_r r^2\left(1+(1+\theta)^2 \right)\right)^2-4\alpha_r(1+\theta)^2 r^2\left(1+ \alpha_r r^2\right) \]
\\
and the roots of this equation will be the following values:
\\
\[ \Delta_{r,1,2} ^2 = \dfrac{1+\theta-\alpha_r r^2\left(1+(1+\theta)^2 \right) \pm \sqrt{D_{\Delta_r^2} }  }{2\alpha_r(1+\theta)^2}.\]
\\
\marginpar{$\textcolor{yellow}{\bullet}$}
From the definition of $ \alpha_r=O_{ r \to 0 }(1) $, that is bounded when $ r $ tends to $ 0 $ and if we consider only a bounded positive number $ \theta $:
\\
\[ \theta>0,\;\;\theta=O_{ r \to 0 }(1), \]
\\
then we have:
\[ D_{\Delta_r^2}= \left( 1+\theta\right)^2+O(r^2).\]
\\
\marginpar{$\textcolor{green}{\bullet}$}
Hence there exists a radius $ 0<R_v \leqslant R_U$ such that for any $ 0 < r \leqslant R_v$:
\[ D_{\Delta_r^2}>0, \]
\\
i.e., for any $ \theta>0,\;\theta=O_{ r \to 0 }(1) $ and any $ 0 < r \leqslant R_v$, there exist such points $ a_r $ and $ v_r $ that the equality \eqref{eq_Eqviv_ar_vr} is satisfied.
\\
\\
\marginpar{$\textcolor{yellow}{\bullet}$}
Let's choose the second root as the solution to the equation \eqref{eq_Ur_t_Dr_ar_fin}:
\[ \Delta_r^2= \dfrac{1+\theta-\alpha_r r^2\left(1+(1+\theta)^2 \right) - \sqrt{D_{\Delta_r^2} }  }{2\alpha_r(1+\theta)^2}.\]
\\
\marginpar{$\textcolor{yellow}{\bullet}$}
And there exists a radius $ 0<R_w \leqslant R_v $ such that for any $ 0 < r \leqslant R_w$ a multiplier\\
\[ 1+\theta-\alpha_r r^2\left(1+(1+\theta)^2 \right) + \sqrt{D_{\Delta_r^2} } \ne 0 \]
\\
and:
\begin{eqnarray} 
&\Delta_r^2= \dfrac{\left( 1+\theta-\alpha_r r^2\left(1+(1+\theta)^2 \right)\right)^2 - D_{\Delta_r^2} }{2\alpha_r(1+\theta)^2\left(1+\theta-\alpha_r r^2\left(1+(1+\theta)^2 \right) + \sqrt{D_{\Delta_r^2} }  \right) }=\nonumber\\
&=\dfrac{4\alpha_r(1+\theta)^2 r^2\left(1+ \alpha_r r^2\right) }{2\alpha_r(1+\theta)^2\left(1+\theta-\alpha_r r^2\left(1+(1+\theta)^2 \right) + \sqrt{D_{\Delta_r^2} }  \right) }=\nonumber\\
&=\dfrac{2r^2\left(1+ \alpha_r r^2\right) }{1+\theta-\alpha_r r^2\left(1+(1+\theta)^2 \right) + \sqrt{D_{\Delta_r^2} } }.\nonumber
\end{eqnarray}
\\
\marginpar{$\textcolor{yellow}{\bullet}$}
And hence at $ r \to 0 $ we have:
\\
\[ \Delta_r^2= \dfrac{2r^2\left(1+ \alpha_r r^2\right) }{1+\theta+O(r^2) + 1+\theta+O(r^2) }=\dfrac{2r^2\left(1+ \alpha_r r^2\right) }{2\left( 1+\theta\right) +O(r^2)}=\dfrac{r^2}{ 1+\theta}+O(r^4).\]
\\
\marginpar{$\textcolor{yellow}{\bullet}$}
Total at $ r \to 0 $:
\begin{align} 
\Delta_r=\dfrac{r}{ \sqrt{1+\theta}}+O(r^4).\nonumber
\end{align}
\\
\marginpar{$\textcolor{yellow}{\bullet}$}
Let's take $\theta $ as the value:
\[ \theta = r^3. \]
\\
\marginpar{$\textcolor{green}{\bullet}$}
Then at $ r \to 0 $:
\begin{align} \label{eq_Delta_r_fin}
\Delta_r=r+O(r^4).
\end{align}
\\
\marginpar{$\textcolor{yellow}{\bullet}$}
And there exists a radius $ 0<R_0 \leqslant R_w $ such that for any $ 0 < r \leqslant R_0$ the expression from the right-hand side of the equality \eqref{eq_Delta_r_fin} will be positive and hence the number $ \theta \Delta_r = r^3\Delta_r$ will be positive, i.e., the interval between points $ a_r $ and $ v_r $ will be nonempty. 
\\
\\
\marginpar{$\textcolor{yellow}{\bullet}$}
Hence from the same mean value theorem and the fulfillment of the equality \eqref{eq_Eqviv_ar_vr}, there must exist a point on the non-zero interval $ (a_r,\; v_r) $:
\[ z_r \stackrel{\scriptscriptstyle\mathrm{def}}{=}x_{z_r}+i\left( t_q+r\right),\;\; x_{a_r} < x_{z_r}< x_{v_r}, \] 
\\
such as this:
\begin{align} \label{eq_diff_zero}
\dfrac{\partial}{\partial x} \left. Re \dfrac{{\zeta }^{\prime }\left(x+i(t_q+r)\right)}{\zeta \left(x+i(t_q+r)\right)}\right|_{ x=x_{z_r} } = 0.
\end{align}
\\
\marginpar{$\textcolor{green}{\bullet}$}
The equality \eqref{eq_diff_zero} with respect to \eqref{eq_Def_Tqs} can be written as follows:
\[ \dfrac{\partial}{\partial x} Re T_q(x+i(t_q+r))\left|_{ x=x_{z_r} }\right.-k(q) \dfrac{(x_{z_r}-\sigma_q)^2-r^2}{\left( (x_{z_r}-\sigma_q)^2+r^2\right)^2 }=0.\]
\\
\marginpar{$\textcolor{green}{\bullet}$}
And then at $ r \to 0 $, given \eqref{eq_Delta_r_fin}:
\begin{eqnarray}
&x_{a_r}=\sigma_q+\Delta_r=\sigma_q+r+O(r^4);\nonumber\\
&x_{v_r}=\sigma_q+\Delta_r+r^3\Delta_r=\sigma_q+r+O(r^4).\nonumber
\end{eqnarray}
\marginpar{$\textcolor{yellow}{\bullet}$}
hence:
\[ x_{z_r}=\sigma_q+r+O(r^4). \]
\\
\marginpar{$\textcolor{yellow}{\bullet}$}
And the equality \eqref{eq_diff_zero} will look like this when $ r \to 0 $:
\begin{eqnarray}
&\dfrac{\partial}{\partial x} Re T_q(x+i(t_q+r))\left|_{ x=x_{z_r} }\right.=\nonumber\\
&\nonumber\\
&=k(q) \dfrac{(r+O(r^4))^2-r^2}{\left( (r+O(r^4))^2+r^2\right)^2 }= \dfrac{O(r^5)}{\left( 2r^2+O(r^5)\right)^2 }=\dfrac{O(r^5)}{4r^4+O(r^7) }=O(r).\nonumber
\end{eqnarray}
\\
\marginpar{$\textcolor{yellow}{\bullet}$}
Then:
\begin{align}
\dfrac{\partial}{\partial x} Re T_q(x+it_q)\left|_{ x=\sigma_q }\right.=\lim_{ r \to 0}\dfrac{\partial}{\partial x} Re T_q(x+i(t_q+r))\left|_{ x=x_{z_r} }\right.=\lim_{ r \to 0}O(r)=0.\nonumber
\end{align}
\\
\marginpar{$\textcolor{yellow}{\bullet}$}
Which contradicts \eqref{eq_NotZero_q}.
\\
\\
\marginpar{$\textcolor{yellow}{\bullet}$}
Hence the assumption \eqref{eq_NotZero_q} is incorrect and as a result:
\\
\begin{align} \label{eq_diff_zero_fin1}
\dfrac{\partial}{\partial x} Re T_q(x+it_q)\left|_{ x=\sigma_q }\right.=0.
\end{align}
\\
\marginpar{$\textcolor{yellow}{\bullet}$}
Let's denote:
\begin{align} 
T_{1-q}(s)\stackrel{\scriptscriptstyle\mathrm{def}}{=}\dfrac{{\zeta }^{\prime }\left(1-s\right)}{\zeta \left(1-s\right)}+\dfrac{k(q)}{s-q} =\dfrac{1}{2}\mathrm{ln}\pi-\dfrac{1}{2}\dfrac{{\Gamma }^{\prime }\left(\dfrac{1-s}{2}\right)}{\Gamma \left(\dfrac{1-s}{2}\right)}-I_{\mathcal{P}\setminus \left\lbrace q\right\rbrace }(s).\nonumber 
\end{align}
\\
\\
\marginpar{$\textcolor{green}{\bullet}$}
The function $ T_{1-q}(s) $ is similarly $ T_q(s) $ analytic and bounded in the entire neighborhood of $Q(R) $.
\\
\\
\marginpar{$\textcolor{green}{\bullet}$}
And the sum of $ T_q(s) $ and $ T_{1-q}(s) $ for any $ Re\left(s\right) >0 $:
\\
\begin{equation}\label{eq_Sum_Tq_1-q}
T_q(s)+T_{1-q}(s)=-\dfrac{1}{2}\dfrac{{\Gamma }^{\prime }\left(\dfrac{s}{2}\right)}{\Gamma \left(\dfrac{s}{2}\right)}-\dfrac{1}{2}\dfrac{{\Gamma }^{\prime }\left(\dfrac{1-s}{2}\right)}{\Gamma \left(\dfrac{1-s}{2}\right)}+\ln\pi. 	
\end{equation}
\\ 
Similarly, all the steps that were done with the functions $ \dfrac{{\zeta }^{\prime }\left(s\right)}{\zeta \left(s\right)} $ and $ T_q(s) $ will be done with the functions $ \dfrac{{\zeta }^{\prime }\left(1-s\right)}{\zeta \left(1-s\right)} $ and $ T_{1-q}(s) $:
\\
\\
\marginpar{$\textcolor{green}{\bullet}$}
Let's define for some $ \Delta_r^*,\;\theta^* \in \mathbb{R}:\;\;  \Delta_r^*>0,\;\;\theta^*>0  $:
\[a_r^* \stackrel{\scriptscriptstyle\mathrm{def}}{=}\sigma_q+\Delta_r^*+i\left( t_q+r\right),  \]
\[v_r^* \stackrel{\scriptscriptstyle\mathrm{def}}{=}\sigma_q+\Delta_r^*+\theta^*\Delta_r^* +i\left( t_q+r\right).  \]
\\
\marginpar{$\textcolor{yellow}{\bullet}$}
Further from the equation:
\\
\[ Re\dfrac{{\zeta }^{\prime }\left(1-a_r^*\right)}{\zeta \left(1-a_r^*\right)}=Re\dfrac{{\zeta }^{\prime }\left(1-v_r^*\right)}{\zeta \left(1-v_r^*\right)}, \]
\\
which through the introduced function $ T_{1-q}(s)$ looks like this: 
\[ ReT_{1-q}(a_r^*)-Re\dfrac{k(q)}{a_r^*-q}= ReT_{1-q}(v_r^*)-Re\dfrac{k(q)}{v_r^*-q}, \]
\\
or:
\[ k(q)\dfrac{\Delta_r^*}{r^2+{\Delta_r^*}^2}-k(q)\dfrac{\Delta_r^*(1+\theta^*)}{r^2+\left( \Delta_r^*(1+\theta^*) \right)^2 }=ReT_{1-q}(a_r^*)-ReT_{1-q}(v_r^*) \]
\\
will follow the expression for the variable $\Delta_r^* $ and the free coefficient $\theta^* $:
\\
\[ \dfrac{1}{r^2+{\Delta_r^*}^2}-\dfrac{1+\theta^*}{r^2+\left( \Delta_r^*(1+\theta^*) \right)^2 }=\theta^*\alpha_r^*, \]
where for some point $ u_r^* \in (a_r^*,\;v_r^*) $:
\\
\[ \alpha_r^*\stackrel{\scriptscriptstyle\mathrm{def}}{=}-\dfrac{1}{k(q)}\dfrac{\partial}{\partial x} Re T_{1-q}(x+i(t_q+r))\left|_{ x=x_{u_r^*} }\right.. \]
\\
\marginpar{$\textcolor{yellow}{\bullet}$}
And it is also shown that when $ r \to 0 $:
\begin{align} 
\Delta_r^*=\dfrac{r}{ \sqrt{1+\theta^*}}+O(r^4).\nonumber
\end{align}
\\
\marginpar{$\textcolor{green}{\bullet}$}
It also takes as its value $\theta^* $:
\[ \theta^* = r^3, \]
and the existence of a solution of the equation under consideration is shown, and at $ r \to 0 $ it is equal:
\begin{align} 
\Delta_r^*=r+O(r^4).\nonumber
\end{align}
\marginpar{$\textcolor{yellow}{\bullet}$}
Next, the existence of the point $ z_r^* \in (a_r^*,\;v_r^*) $ follows from the mean value theorem:
\\
\[ \dfrac{\partial}{\partial x} \left. Re \dfrac{{\zeta }^{\prime }\left(1-x-i(t_q+r)\right)}{\zeta \left(1-x-i(t_q+r)\right)}\right|_{ x=x_{z_r} } = 0, \]
\\
which in turn leads to an expression similar to \eqref{eq_diff_zero_fin1}:
\begin{align} \label{eq_diff_zero_fin2}
\dfrac{\partial}{\partial x} Re T_{1-q}(x+it_q)\left|_{ x=\sigma_q }\right.=0.
\end{align}
\\
\marginpar{$\textcolor{yellow}{\bullet}$}
Adding the left and right parts of the equations \eqref{eq_diff_zero_fin1} and \eqref{eq_diff_zero_fin2}, taking into account \eqref{eq_Sum_Tq_1-q}, we obtain:
\[ \dfrac{\partial}{\partial x}\left. \left( Re T_q(x+it_q) + Re T_{1-q}(x+it_q)\right)\right|_{ x=\sigma_q } =0, \]
\\
or total:
\\
\[\dfrac{\partial}{\partial x}\left.Re\left( -\dfrac{1}{2}\dfrac{{\Gamma }^{\prime }\left(\dfrac{x+it_q}{2}\right)}{\Gamma \left(\dfrac{x+it_q}{2} \right)} -\dfrac{1}{2}\dfrac{{\Gamma }^{\prime }\left(\dfrac{1-x-it_q}{2}\right)}{\Gamma \left(\dfrac{1-x-it_q}{2} \right)}\right)\right|_{x=\sigma_q}=0.\;\;\;\square \]
\\
STATEMENT B
\[ \forall\; 0<\sigma \leqslant \dfrac{1}{2},\;y\geqslant 4: \]
\marginpar{$\textcolor{red}{\bullet}$}
\begin{eqnarray} \label{eq_gamma_1_2}
	&\dfrac{\partial}{\partial x}\left.Re\left( -\dfrac{1}{2}\dfrac{{\Gamma }^{\prime }\left(\dfrac{x+iy}{2}\right)}{\Gamma \left(\dfrac{x+iy}{2} \right)} -\dfrac{1}{2}\dfrac{{\Gamma }^{\prime }\left(\dfrac{1-x-iy}{2}\right)}{\Gamma \left(\dfrac{1-x-iy}{2} \right)}\right)\right|_{x=\sigma}=0 \Leftrightarrow \\
	&\nonumber\\
	&\nonumber\\
	&\Leftrightarrow \sigma =\dfrac{1}{2}.\nonumber
\end{eqnarray}
\\
PROOF:
\\
\\
\marginpar{$\textcolor{green}{\bullet}$}
Given \eqref{eq_PrimeLn}, \eqref{eq_PrimeLnDigamma}, and the Digamma function formula from~\cite[p.259 \S 6.3.16]{Abramowitz-1972}:
\begin{eqnarray} 
	&Re\left(  -\dfrac{1}{2}\dfrac{{\Gamma }^{\prime }\left(\dfrac{x+iy}{2}\right)}{\Gamma \left(\dfrac{x+iy}{2} \right)}-\dfrac{1}{2}\dfrac{{\Gamma }^{\prime }\left(\dfrac{1-x-iy}{2}\right)}{\Gamma \left(\dfrac{1-x-iy}{2} \right)} \right)=\nonumber\\
	&=Re\left(  \dfrac{\gamma}{2}+\dfrac{1}{x+iy}+\sum _{n=1}^{\infty }\left(\dfrac{1}{x+iy+2n}-\dfrac{1}{2n}\right)+\right. \nonumber \\
	&\left. +\dfrac{\gamma}{2}+ \dfrac{1}{1-x-iy}+\sum _{n=1}^{\infty }\left(\dfrac{1}{1-x-iy+2n}-\dfrac{1}{2n}\right)\right)= \nonumber
\end{eqnarray}
\begin{eqnarray} \label{eq_sum_gamma_fin2}
	&=  \gamma +\dfrac{x}{x^2+y^2}+\sum _{n=1}^{\infty }\left(\dfrac{2n+x}{(2n+x)^2+y^2}-\dfrac{1}{2n}\right)+ \nonumber \\
	&+ \dfrac{1-x}{(1-x)^2+y^2} +\sum _{n=1}^{\infty }\left(\dfrac{2n+1-x}{(2n+1-x)^2+y^2}-\dfrac{1}{2n}\right).
\end{eqnarray}
\\
\\
\\
\\
\\
\marginpar{$\textcolor{yellow}{\bullet}$}
From \eqref{eq_sum_gamma_fin2}, the equality \eqref{eq_gamma_1_2} can be viewed as an equation with respect to variable $ x $ with constraint $ 0<x \leqslant \dfrac{1}{2},\;y\geqslant 4 $:
\begin{eqnarray} \label{eq_DivXX_0}
	\sum _{n=0}^{\infty } \left(  \dfrac{(2n+x)^2-y^2}{((2n+x)^2+y^2)^2}-\dfrac{(2n+1-x)^2-y^2}{((2n+1-x)^2+y^2)^2}\right)=0
\end{eqnarray}
\\
\\
\marginpar{$\textcolor{green}{\bullet}$}
In turn:
\begin{eqnarray}
	&\sum _{n=0}^{\infty } \left(  \dfrac{(2n+x)^2-y^2}{((2n+x)^2+y^2)^2}-\dfrac{(2n+1-x)^2-y^2}{((2n+1-x)^2+y^2)^2}\right)=\nonumber\\
	&=\sum _{n=0}^{\infty } \left(  \dfrac{1}{(2n+x)^2+y^2}-\dfrac{1}{(2n+1-x)^2+y^2}\right)-\nonumber\\
	&-2y^2\sum _{n=0}^{\infty } \left(  \dfrac{1}{((2n+x)^2+y^2)^2}-\dfrac{1}{((2n+1-x)^2+y^2)^2}\right)=\nonumber
\end{eqnarray}
\begin{eqnarray} \label{eq_DivXX_line}
	&=\sum _{n=0}^{\infty }   \dfrac{(1-2x)(4n+1)}{((2n+x)^2+y^2)((2n+1-x)^2+y^2)}-\nonumber\\
	&-2y^2\sum _{n=0}^{\infty }   \dfrac{(1-2x)(4n+1)((2n+x)^2+(2n+1-x)^2+2y^2)}{((2n+x)^2+y^2)^2((2n+1-x)^2+y^2)^2}=\nonumber\\
	&=(1-2x) \left( \sum _{n=0}^{\infty }   \dfrac{4n+1}{((2n+x)^2+y^2)((2n+1-x)^2+y^2)}-\right. \nonumber\\
	&\left. -2y^2\sum _{n=0}^{\infty }   \dfrac{(4n+1)((2n+x)^2+(2n+1-x)^2+2y^2)}{((2n+x)^2+y^2)^2((2n+1-x)^2+y^2)^2}\right).
\end{eqnarray}
\\
\\
\marginpar{$\textcolor{green}{\bullet}$}
Let us evaluate the sum in the general brackets of the equality \eqref{eq_DivXX_line}:
\begin{eqnarray}
	&\sum_{n=0}^{\infty }   \dfrac{4n+1}{((2n+x)^2+y^2)((2n+1-x)^2+y^2)}-\nonumber\\
	& -2y^2\sum _{n=0}^{\infty }   \dfrac{(4n+1)((2n+x)^2+(2n+1-x)^2+2y^2)}{((2n+x)^2+y^2)^2((2n+1-x)^2+y^2)^2}.\nonumber
\end{eqnarray}
\\
\\
\\
\marginpar{$\textcolor{green}{\bullet}$}
From~\cite[p.259]{Abramowitz-1972}, \cite[\S~6.495]{Adams-1922} :
\[ \sum _{n=1}^{\infty } \dfrac{1}{(2n-1)^2+y^2}=\dfrac{\pi}{4y}\tanh \dfrac{\pi y}{2}, \]
\[ \sum _{n=1}^{\infty } \dfrac{1}{(2n)^2+y^2}=-\dfrac{1}{2y^2}+\dfrac{\pi}{4y}\coth \dfrac{\pi y}{2}.\]
And then the first summand in the sum in question:
\begin{eqnarray}
	&\sum_{n=0}^{\infty }   \dfrac{4n+1}{((2n+x)^2+y^2)((2n+1-x)^2+y^2)}<\nonumber\\
	&<\dfrac{1}{(x^2+y^2)((1-x)^2+y^2)}+\sum_{n=1}^{\infty } \dfrac{4n+1}{((2n-1)^2+y^2)((2n)^2+y^2)}<\nonumber\\
	&<\dfrac{1}{y^4}+\sum_{n=1}^{\infty }\left( \dfrac{1}{(2n-1)^2+y^2}-\dfrac{1}{(2n)^2+y^2}\right)+\nonumber\\
	&+\sum_{n=1}^{\infty }\dfrac{2}{((2n-1)^2+y^2)^2}.\nonumber
\end{eqnarray}
\marginpar{$\textcolor{green}{\bullet}$}
Here:
\[\sum _{n=1}^{\infty }\left( \dfrac{1}{(2n-1)^2+y^2}-\dfrac{1}{(2n)^2+y^2}\right)=\dfrac{\pi}{4y}\tanh \dfrac{\pi y}{2}-\dfrac{\pi}{4y} \coth \dfrac{\pi y}{2} +\dfrac{1}{2y^2}, \]
\begin{eqnarray}
	&\sum _{n=1}^{\infty }\dfrac{2}{((2n-1)^2+y^2)^2}=-\dfrac{1}{y}\dfrac{d}{d y}\left( \sum _{n=1}^{\infty } \dfrac{1}{(2n-1)^2+y^2}\right)=\nonumber\\
	&=-\dfrac{1}{y}\dfrac{d}{d y}\left(\dfrac{\pi}{4y}\tanh \dfrac{\pi y}{2} \right)= \dfrac{\pi}{4y^3}\tanh \dfrac{\pi y}{2}-\dfrac{\pi^2}{8y^2}\dfrac{1}{\cosh^2 \dfrac{\pi y}{2}} \nonumber
\end{eqnarray}
I.e.
\begin{eqnarray} 
	&\sum_{n=0}^{\infty }   \dfrac{4n+1}{((2n+x)^2+y^2)((2n+1-x)^2+y^2)}<\nonumber\\
	&<\dfrac{1}{2y^2}+\dfrac{\pi}{4y^3}\tanh \dfrac{\pi y}{2}+\dfrac{1}{y^4}-\nonumber\\
	&-\dfrac{\pi}{4y}\left( \coth \dfrac{\pi y}{2}-\tanh \dfrac{\pi y}{2}\right)-\dfrac{\pi^2}{8y^2}\dfrac{1}{\cosh^2 \dfrac{\pi y}{2}}. \nonumber
\end{eqnarray}
\marginpar{$\textcolor{green}{\bullet}$}
Second summation:
\begin{eqnarray} 
	&\sum _{n=0}^{\infty }   \dfrac{(4n+1)((2n+x)^2+(2n+1-x)^2+2y^2)}{((2n+x)^2+y^2)^2((2n+1-x)^2+y^2)^2}=\nonumber\\
	& =\sum _{n=1}^{\infty }   \dfrac{4n-3}{((2n-2+x)^2+y^2)((2n-1-x)^2+y^2)}*\nonumber\\
	&*\left(\dfrac{1}{(2n-2+x)^2+y^2}+\dfrac{1}{(2n-1-x)^2+y^2} \right)>\nonumber\\
	&>\sum _{n=1}^{\infty }   \dfrac{4n-1}{((2n-1)^2+y^2)((2n)^2+y^2)}\left(\dfrac{1}{(2n-1)^2+y^2}+\dfrac{1}{(2n)^2+y^2} \right)-\nonumber\\
	&-\sum _{n=1}^{\infty }   \dfrac{2}{((2n-1)^2+y^2)((2n)^2+y^2)}\left(\dfrac{1}{(2n-1)^2+y^2}+\dfrac{1}{(2n)^2+y^2} \right)>\nonumber\\
	&>\sum _{n=1}^{\infty }  \left(\dfrac{1}{((2n-1)^2+y^2)^2}-\dfrac{1}{((2n)^2+y^2)^2} \right) -\nonumber\\
	&-\sum _{n=1}^{\infty } \dfrac{4}{((2n-1)^2+y^2)^3}.\nonumber
\end{eqnarray} 
\marginpar{$\textcolor{green}{\bullet}$}
Here:
\begin{eqnarray}
	&\sum _{n=1}^{\infty }  \left(\dfrac{1}{((2n-1)^2+y^2)^2}-\dfrac{1}{((2n)^2+y^2)^2} \right)=\nonumber\\
	&=-\dfrac{1}{2y}\dfrac{d}{d y}\left( \sum _{n=1}^{\infty } \dfrac{1}{(2n-1)^2+y^2}-\dfrac{1}{(2n)^2+y^2}\right)=\nonumber\\
	&=-\dfrac{1}{2y}\dfrac{d}{d y}\left(\dfrac{\pi}{4y}\tanh \dfrac{\pi y}{2} +\dfrac{1}{2y^2}-\dfrac{\pi}{4y}\coth \dfrac{\pi y}{2}\right)=\nonumber\\
	&= \dfrac{\pi}{8y^3}\tanh \dfrac{\pi y}{2}-\dfrac{\pi^2}{16y^2}\dfrac{1}{\cosh^2 \dfrac{\pi y}{2}} +\dfrac{1}{2y^4}-\dfrac{\pi}{8y^3}\coth \dfrac{\pi y}{2}-\dfrac{\pi^2}{16y^2}\dfrac{1}{\sinh^2 \dfrac{\pi y}{2}}.\nonumber
\end{eqnarray}
\marginpar{$\textcolor{green}{\bullet}$}
And
\begin{eqnarray}
	&\sum _{n=1}^{\infty } \dfrac{4}{((2n-1)^2+y^2)^3}=\nonumber\\
	&=-\dfrac{1}{2y}\dfrac{d}{dy}\left(-\dfrac{1}{y}\dfrac{d}{dy}\left( \sum _{n=1}^{\infty } \dfrac{1}{(2n-1)^2+y^2} \right)  \right)=\nonumber
\end{eqnarray}
\begin{eqnarray}
	&=-\dfrac{1}{2y}\dfrac{d}{dy}\left(-\dfrac{1}{y}\dfrac{d}{dy}\left( \dfrac{\pi}{4y}\tanh \dfrac{\pi y}{2}\right)  \right)=\nonumber\\
	&=-\dfrac{1}{2y}\dfrac{d}{dy}\left( \dfrac{\pi}{4y^3}\tanh \dfrac{\pi y}{2}-\dfrac{\pi^2}{8y^2}\dfrac{1}{\cosh^2 \dfrac{\pi y}{2}}\right) =\nonumber\\
	&=\dfrac{3\pi}{8y^5}\tanh \dfrac{\pi y}{2}-\dfrac{\pi^2}{16y^4}\dfrac{1}{\cosh^2 \dfrac{\pi y}{2}}-\dfrac{\pi^2}{8y^4}\dfrac{1}{\cosh^2 \dfrac{\pi y}{2}}-\dfrac{\pi^3}{16y^3}\dfrac{\tanh \dfrac{\pi y}{2}}{\cosh^2 \dfrac{\pi y}{2}}=\nonumber\\
	&=-\dfrac{\pi^2}{8y^3 \cosh^2 \dfrac{\pi y}{2}}\left( \dfrac{3}{2y}+\dfrac{\pi}{2} \tanh \dfrac{\pi y}{2}  \right) +\dfrac{3\pi}{8y^5}\tanh \dfrac{\pi y}{2}. \nonumber
\end{eqnarray}
\\
\\
\marginpar{$\textcolor{yellow}{\bullet}$}
Hence:
\begin{eqnarray}
	&\sum_{n=0}^{\infty }   \dfrac{4n+1}{((2n+x)^2+y^2)((2n+1-x)^2+y^2)}-\nonumber\\
	& -2y^2\sum _{n=0}^{\infty }   \dfrac{(4n+1)((2n+x)^2+(2n+1-x)^2+2y^2)}{((2n+x)^2+y^2)^2((2n+1-x)^2+y^2)^2}<\nonumber\\
	&<\dfrac{1}{2y^2}+\dfrac{\pi}{4y^3}\tanh \dfrac{\pi y}{2}+\dfrac{1}{y^4}-\nonumber\\
	&-\dfrac{\pi}{4y}\left( \coth \dfrac{\pi y}{2}-\tanh \dfrac{\pi y}{2}\right)-\dfrac{\pi^2}{8y^2}\dfrac{1}{\cosh^2 \dfrac{\pi y}{2}}-\nonumber\\
	& -2y^2\left(\dfrac{\pi}{8y^3}\tanh \dfrac{\pi y}{2}-\dfrac{\pi^2}{16y^2}\dfrac{1}{\cosh^2 \dfrac{\pi y}{2}} +\dfrac{1}{2y^4}-\right.  \nonumber\\
	&\left.  -\dfrac{\pi}{8y^3}\coth \dfrac{\pi y}{2}-\dfrac{\pi^2}{16y^2}\dfrac{1}{\sinh^2 \dfrac{\pi y}{2}} \right)-\nonumber\\
	&-\left(-\dfrac{\pi^2}{8y^3 \cosh^2 \dfrac{\pi y}{2}}\left( \dfrac{3}{2y}+\dfrac{\pi}{2} \tanh \dfrac{\pi y}{2}  \right) +\dfrac{3\pi}{8y^5}\tanh \dfrac{\pi y}{2} \right)=\nonumber
\end{eqnarray}
\begin{eqnarray} \label{eq_First_delta}
	&=\dfrac{1}{y^2}\left( -\dfrac{1}{2}-\dfrac{\pi^2}{8}\dfrac{1}{\cosh^2 \dfrac{\pi y}{2}}-\dfrac{3\pi}{8y^3}\tanh \dfrac{\pi y}{2}+\right. \nonumber\\
	&\left. +\dfrac{\pi^2}{8}\dfrac{y^2}{\cosh^2 \dfrac{\pi y}{2}}+\dfrac{\pi^2}{8}\dfrac{y^2}{\sinh^2 \dfrac{\pi y}{2}}+\right. \nonumber\\
	&\left. +\dfrac{\pi}{4y}\tanh \dfrac{\pi y}{2}+\dfrac{1}{y^2}+\dfrac{\pi^2}{8y \cosh^2 \dfrac{\pi y}{2}}\left( \dfrac{3}{2y}+\dfrac{\pi}{2} \tanh \dfrac{\pi y}{2}  \right)\right).
\end{eqnarray}
\\
\\
Consider the positive summands inside the common bracket of the right-hand side of the inequality \eqref{eq_First_delta} at $ y \geqslant 4 $:
\\
\\
\marginpar{$\textcolor{green}{\bullet}$}
Derivative:
\begin{eqnarray}
	\left( \dfrac{y^2}{\cosh^2 \dfrac{\pi y}{2}}\right)^{\prime}=y\dfrac{2 \cosh \dfrac{\pi y}{2} - \pi y \sinh \dfrac{\pi y}{2} }{\cosh^3 \dfrac{\pi y}{2}}<0,\nonumber
\end{eqnarray}
since for $ \forall\; y \geqslant 4 $
\[ \dfrac{2}{\pi} \coth \dfrac{\pi y}{2}<y. \]
\marginpar{$\textcolor{green}{\bullet}$}
Similarly, the derivative:
\begin{eqnarray}
	\left( \dfrac{y^2}{\sinh^2 \dfrac{\pi y}{2}}\right)^{\prime}=y\dfrac{2 \sinh \dfrac{\pi y}{2} - \pi y \cosh \dfrac{\pi y}{2} }{\sinh^3 \dfrac{\pi y}{2}}<0,\nonumber
\end{eqnarray}
since for $ \forall\; y \geqslant 4 $
\[ \dfrac{2}{\pi} \tanh \dfrac{\pi y}{2}<y. \]
\marginpar{$\textcolor{green}{\bullet}$}
Hence $ \forall\; y \geqslant 4 $:
\begin{eqnarray}
	\dfrac{\pi^2}{8}\dfrac{y^2}{\cosh^2 \dfrac{\pi y}{2}} \leqslant \dfrac{\pi^2}{8}\dfrac{16}{\cosh^2 2\pi}<0,0002754,\nonumber\\
	\dfrac{\pi^2}{8}\dfrac{y^2}{\sinh^2 \dfrac{\pi y}{2}} \leqslant \dfrac{\pi^2}{8}\dfrac{16}{\sinh^2 2\pi}<0,0002754.\nonumber
\end{eqnarray}
\marginpar{$\textcolor{green}{\bullet}$}
Further $ \forall\; y \geqslant 4 $:
\[ \dfrac{\pi}{4y}\tanh \dfrac{\pi y}{2} < \dfrac{\pi}{16}< 0,1963496, \]
\[ \dfrac{1}{y^2} \leqslant 0,0625, \]
\[ \dfrac{\pi^2}{8y \cosh^2 \dfrac{\pi y}{2}}\left( \dfrac{3}{2y}+\dfrac{\pi}{2} \tanh \dfrac{\pi y}{2}  \right) < \dfrac{\pi^2}{32 \cosh^2 2\pi}\left( \dfrac{3}{8}+\dfrac{\pi}{2} \right)<0,0000084. \]
\\
\\
\\
\\
\marginpar{$\textcolor{yellow}{\bullet}$}
Hence $ \forall\; y \geqslant 4 $ the total sum of positive terms in the total bracket does not exceed $ \dfrac{1}{2} $:
\begin{eqnarray} 
	& \dfrac{\pi^2}{8}\dfrac{y^2}{\cosh^2 \dfrac{\pi y}{2}}+\dfrac{\pi^2}{8}\dfrac{y^2}{\sinh^2 \dfrac{\pi y}{2}}+ \nonumber\\
	& +\dfrac{\pi}{4y}\tanh \dfrac{\pi y}{2}+\dfrac{1}{y^2}+\dfrac{\pi^2}{8y \cosh^2 \dfrac{\pi y}{2}}\left( \dfrac{3}{2y}+\dfrac{\pi}{2} \tanh \dfrac{\pi y}{2}  \right)<0,2594088.\nonumber
\end{eqnarray}
\\
\\
This means that when $ \forall\; y \geqslant 4,\;\;0<x \leqslant \dfrac{1}{2} $ the second multiplier of the right-hand side of the equality \eqref{eq_DivXX_line} does not turn to $ 0 $,
\\
hence from \eqref{eq_DivXX_0} and \eqref{eq_DivXX_line}:
\[ x=\dfrac{1}{2}. \]
\\
\marginpar{$\textcolor{green}{\bullet}$}
On the flip side, the truth of the statement is obvious.
\\
\\
$ \square $
\\
\\
\\
\\
\pagebreak
\\
For $\forall\; \rho=\sigma+it \in \mathcal{P} $, we evaluate the minimum value of $ \left| t\right| $.
\\
\\
\marginpar{$\textcolor{yellow}{\bullet}$}
Let's denote by $ t_1 \stackrel{\scriptscriptstyle\mathrm{def}}{=} \min_{\rho \in\mathcal{P}} \left| Im (\rho )\right| $.
\\
\\
\marginpar{$\textcolor{green}{\bullet}$}
For $ \forall\; \rho=\sigma+it \in \mathcal{P}_1  \cup  \mathcal{P}_3:  $
\\
\begin{eqnarray} 
&\dfrac{1}{\rho}+\dfrac{1}{\bar{\rho}}+\dfrac{1}{1-\rho}+\dfrac{1}{1-\bar{\rho}}=\nonumber\\
&=\dfrac{\sigma}{\sigma^2+t^2}+\dfrac{\sigma}{\sigma^2+t^2}+\dfrac{1-\sigma}{(1-\sigma)^2+t^2}+\dfrac{1-\sigma}{(1-\sigma)^2+t^2}=\nonumber\\
&=\dfrac{2\sigma}{\sigma^2+t^2}+\dfrac{2(1-\sigma)}{(1-\sigma)^2+t^2}>\dfrac{2\sigma}{1+t^2}+\dfrac{2(1-\sigma)}{1+t^2}=\dfrac{2}{1+t^2}.\nonumber
\end{eqnarray}
\\
\\
\marginpar{$\textcolor{green}{\bullet}$}
For $ \forall\; \rho=\dfrac{1}{2}+it \in \mathcal{P}_2:  $
\\
\begin{eqnarray} 
&\dfrac{1}{\rho}+\dfrac{1}{1-\rho}=\nonumber\\
&=\dfrac{\dfrac{1}{2}}{\dfrac{1}{4}+t^2}+\dfrac{\dfrac{1}{2}}{\dfrac{1}{4}+t^2}=
\dfrac{1}{\dfrac{1}{4}+t^2}.\nonumber
\end{eqnarray}
\\
\\
Then taking into account \eqref{eq_RHO} either:
\[ \dfrac{2}{1+t_1^2}<\sum _{\rho \in \mathcal{P}}\dfrac{1}{\rho }<0,0230958, \]
or:
\[ \dfrac{1}{\dfrac{1}{4}+t_1^2}<\sum _{\rho \in \mathcal{P}}\dfrac{1}{\rho }<0,0230958. \]
In the first case:
\begin{equation} 
t_1>9,2518014,\nonumber
\end{equation}
in the second:
\begin{equation} 
t_1>6,5610909.\nonumber
\end{equation}
\\
\marginpar{$\textcolor{green}{\bullet}$}
Therefore:
\begin{equation} \label{eq_t1_rez}
t_1>6,5610909.\nonumber
\end{equation}
\\
\\
\marginpar{$\textcolor{yellow}{\bullet}$}
So when we take an arbitrary nontrivial root $ q=\sigma_q+i t_q \in \mathcal{P}_1 \cup \mathcal{P}_2 $ based on Statements A and B given that $ \left| t_q\right| \geqslant t_1 > 4 $, we conclude that:
\[ \sigma_q=\dfrac{1}{2}, \]
i.e.:
\[ \mathcal{P}_1=\varnothing. \] 
\\
\\
\marginpar{$\textcolor{yellow}{\bullet}$}
And given that $\left\||\mathcal{P}3 \right\| =\left\|\mathcal{P}_1 \right\|=0$, we have:
\\
\[ \mathcal{P}_3=\mathcal{P}_1=\varnothing \]
and
\\
\[\mathcal{P}=\mathcal{P}_2,\]
\\
which proves the main statement and conjecture put forward by Bernhard Riemann about the location of the real parts of nontrivial zeros of the Zeta-function.
\\
\\
$ \square $
\\
\pagebreak
\renewcommand{\refname}{REFERENCES:}

\end{document}